\documentclass[a4paper,3p]{elsarticle}
\usepackage[OT4]{fontenc}
\usepackage{amssymb}
\usepackage{amsmath}
\usepackage[all,2cell]{xy}
\UseAllTwocells

\newcommand{\cat}[1]{\ensm{\textup{\bf #1}}}
\newcommand{\loc}[1]{\ensm{\opm{\cat{#1}}}}
\newcommand{\topt}[1]{\ensm{\cat{Af\,Spc}(#1)}}
\newcommand{\set}{\ensm{\cat{Set}}}
\newcommand{\algo}{\ensm{\cat{Alg}(\Omega)}}
\newcommand{\algop}[1]{\ensm{\cat{Alg}(#1)}}
\newcommand{\topsyst}[1]{\ensm{\cat{Af\,Sys}(#1)}}
\newcommand{\spt}[1]{\ensm{\cat{Af\,Sys}_s(#1)}}
\newcommand{\semilatc}[1]{\ensm{\cat{CSLat}(#1)}}
\newcommand{\btop}[1]{\ensm{{#1}\text{-}\cat{Top}}}
\newcommand{\btopsys}[1]{\ensm{{#1}\text{-}\cat{TopSys}}}
\newcommand{\afinst}[1]{\ensm{\cat{AfInst}(#1)}}
\newcommand{\safinst}[1]{\ensm{\cat{SAfInst}(#1)}}
\newcommand{\lafinst}[1]{\ensm{\cat{LAfInst}(#1)}}
\newcommand{\syst}[1]{\ensm{\cat{GAfSys}(#1)}}

\newcommand{\arw}[1]{\ensm{\xrightarrow{#1}}}
\newcommand{\incl}[3]{\ensm{\xymatrix{{#1\,}\ar@{^{(}->}[r]^-{#2} & {#3}}}}
\newcommand{\cell}[4]{\ensm{\xymatrix{{#1} \ar@<0.5ex>[r]^-{#2} \ar@<-0.5ex>[r]_-{#3} & {#4}}}}

\newcommand{\sqed}{\hfill{\vrule width 3pt height 3pt depth 0pt}}

\newcommand{\ensm}[1]{\ensuremath{#1}}
\newcommand{\opm}[1]{\ensm{{#1}^{op}}}
\newcommand{\mcal}[1]{\ensm{\mathcal{#1}}}
\newcommand{\msf}[1]{\ensm{\mathsf{#1}}}
\newcommand{\seq}{\ensm{\subseteq}}
\newcommand{\ovr}[1]{\ensuremath{\overline{#1}}}
\newcommand{\prm}[1]{\ensm{{#1}^{\prime}}}
\newcommand{\obc}[1]{\ensm{Ob(#1)}}
\newcommand{\mbb}[1]{\ensm{\mathbb{#1}}}
\newcommand{\leqs}{\ensm{\leqslant}}

\DeclareMathOperator{\pt}{pt}
\DeclareMathOperator{\omg}{\mathit{\Omega}}

\journal{Fuzzy Sets and Systems}

\begin{document}

\begin{frontmatter}

\title{Topological systems as a framework for institutions}

\author[Jeffrey]{Jeffrey T. Denniston}
\address[Jeffrey]{Department of Mathematical Sciences, Kent State University\\
                  Kent, Ohio, USA 44242}
\ead{jdennist@kent.edu}

\author[Austin]{Austin Melton}
\address[Austin]{Departments of Computer Science and Mathematical Sciences, Kent State University\\
                  Kent, Ohio, USA 44242}
\ead{amelton@kent.edu}

\author[Stephen]{Stephen E. Rodabaugh}
\address[Stephen]{College of Science, Technology, Engineering, Mathematics (STEM), Youngstown State University\\
                  Youngstown, Ohio, USA 44555-3347\\}
\ead{serodabaugh@ysu.edu}

\author[Sergey]{Sergey A. Solovyov\tnoteref{GACR}}
\address[Sergey]{Institute of Mathematics, Faculty of Mechanical Engineering, Brno University of Technology, Technicka 2896/2,
616 69 Brno, Czech Republic}
\ead{solovjovs@fme.vutbr.cz}
\tnotetext[GACR]{The author gratefully acknowledges the support of Grant Agency of Czech Republic (GA\v{C}R) and
Austrian Science Fund (FWF) within bilateral project No. I 1923-N25 ``New Perspectives on Residuated Posets".}

\begin{abstract}
 Recently, J.~T.~Denniston, A.~Melton, and S.~E.~Rodabaugh introduced a lattice-valued analogue of the concept of institution of J.~A.~Goguen and R.~M.~Burstall, comparing it, moreover, with the (lattice-valued version of the) notion of topological system of S.~Vickers. In this paper, we show that a suitable generalization of topological systems provides a convenient framework for doing certain kinds of (lattice-valued) institutions.
\end{abstract}

\begin{keyword} adjoint situation \sep affine theory \sep comma category \sep elementary institution \sep localification and spatialization procedure \sep topological institution \sep topological space \sep topological system \sep variety of algebras

 \MSC[2010] 18A25 \sep 18B15 \sep 18B30 \sep 18B99 \sep 18C10
\end{keyword}

\end{frontmatter}


\newtheorem{thm}{Theorem}
\newtheorem{prop}[thm]{Proposition}
\newtheorem{cor}[thm]{Corollary}
\newtheorem{lem}[thm]{Lemma}
\newdefinition{defn}[thm]{Definition}
\newdefinition{exmp}[thm]{Example}
\newdefinition{rem}[thm]{Remark}
\newdefinition{prob}[thm]{Problem}
\newproof{pf}{Proof}

\section{Introduction}


\par
There exists a convenient approach to logical systems in computer science, which is based in the concept of \emph{institution} of J.~A.~Goguen and R.~M.~Burstall~\cite{Goguen1984}. An institution comprises a category of (abstract) signatures, where every signature has its associated sentences, models, and a relationship of satisfaction; this relationship is invariant (in a certain sense) under change of signature. The slogan, therefore, is ``truth is invariant under change of notation". Examples of institutions include unsorted universal algebra, many-sorted algebra, order-sorted algebra, several variants of first-order logic, partial algebra (see, e.g.,~\cite{Goguen2002}). More examples can be found in~\cite[Subsection~3.2]{Diaconescu2008}. A number of authors (including the initiators themselves) have proposed generalizations of institutions in various forms as well as advanced their theories~\cite{Goguen1986,Goguen1992,Mayoh1985,Mossakowski1996,Salibra1993}. Moreover, some authors used a purely category-theoretic approach to institutions (see, e.g.,~\cite{Diaconescu2002}).

\par
There exists the concept of \emph{topological system} of S.~Vickers~\cite{Vickers1989}, which is based in the ideas of geometric
logic~\cite{Vickers1993a} and intended to provide a common setting for both topological spaces (point-set topology) and their underlying algebraic structures\textrm{---}locales (point-free topology). In particular, S.~Vickers presented system spatialization and localification procedures, which created ways to move back and forth between the categories of topological spaces (resp., locales) and topological systems. Recently, the latter concept has gained in interest in connection with lattice-valued topology. In particular,~\cite{Denniston2009,Denniston2012} introduced and studied the notion of lattice-valued topological system;~\cite{Guido2010} discovered a convenient relationship between crisp and lattice-valued topology, based in topological systems; and~\cite{Solovjovs2009,Solovyov2012} studied a lattice-valued analogue of the above-mentioned system spatialization procedure.

\par
In an attempt to find possible relationships between institutions and topological systems, at the 35th Linz Seminar on Fuzzy Set Theory, J.~T.~Denniston, A.~Melton, and S.~E.~Rodabaugh~\cite{Denniston2014} presented a lattice-valued analogue of institutions, and showed that (lattice-valued) topological systems provide a particular instance of the latter. Moreover,~\cite{Sernadas1995} introduced (crisp) \emph{topological institutions}, based in topological systems, the slogan being that ``the central concept is the theory, not the formula". To continue this line of study, several authors considered some other institutional modifications (e.g., probability institutions, quantum institutions, etc.)~\cite{Baltazar2006,Caleiro2006}, motivated by the ideas of quantum logic (in connection with quantum physics).

\par
We notice that some researchers prefer the reverse of the above slogan. While the theory plays an important role in building logic, the terms should be constructed first (to allow variable substitution), and after that sentences should be constructed on terms to get a sentence functor. The construction of terms through the term monad plays a key role in allowing variable substitution and variable assignment. Having only abstract categories for signatures hides the use of variables and the difference between terms and sentences. Term and sentence construction are two separate processes, which should be revealed together with the process of variable assignment. The readers with the same point of view could look into~\cite{Eklund2014,Eklund2002,Eklund2004,Eklund2007,Eklund2012} for a particular fuzzy approach to terms and their respective monad. We notice, however, that although the main point of institution theory is exactly to liberate the logic study from explicit variables and substitutions, when needed institution theory has its well established approach to these~\cite{Diaconescu2008}.

\par
The main purpose of this paper is to show that a suitably generalized concept of topological system makes a setting for a particular type of (lattice-valued) institutions, namely, elementary institutions~\cite{Sernadas1994,Sernadas1995}. In so doing, we aim at providing a convenient framework for building up the theory of lattice-valued institutions. More precisely, there already exists a well-developed theory of lattice-valued (also many-valued or fuzzy) logic, which has been given a coherent statement by P.~H\'ajek in~\cite{Hajek1998}, and which by now has much diversified w.r.t. the algebraic structures (which often constitute a variety) over which the respective fuzzification is done. The concept of institution, however, being a significant part of the crisp logical developments (see, e.g.,~\cite{Diaconescu2008}), has been fuzzified just recently in~\cite{Denniston2014}. With this article, we are going to extend further this fuzzification in a way, which could encompass various lattice-valued frameworks. We achieve this goal with a modification of the affine context of Y.~Diers~\cite{Diers1996,Diers1999,Diers2002}, which is based in an arbitrary variety of algebras, thereby providing a unifying setting for many possible fuzzifications of institutions (which are to come), each of them based in the favourite variety of it's authors (e.g., a variety of residuated lattices). The main advantage of such a unifying setting is the fact that every statement, which is proved in the affine framework (namely, for all varieties) will be valid for each particular fuzzification (namely, for each particular variety).


\section{Affine systems and their related tools}


\par
This section reviews the notions of affine systems and spaces, as well as their related spatialization and localification procedures (for details, see~\cite{Solovyovi,Solovyov2012,Solovyov2013a}). Since the localification procedure has not yet appeared in the literature in full detail, we provide a more thorough description. We conclude this section with the approach to topological systems, which is motivated by algebraic theories of~F.~W.~Lawvere~\cite{Lawvere1963}.

\par
A particular remark is helpful w.r.t. the system terminology of this paper. Following the notion of lattice-valued topological system of~\cite{Denniston2009,Denniston2012}, in~\cite{Solovyov2012} is provided a more general concept under the name of variety-based topological system, which eventually gave rise to the notion and theory of categorically-algebraic topology~\cite{Solovjovs2010c}. It was subsequently discovered that the latter concept had already been introduced by Y.~Diers~\cite{Diers1996,Diers1999,Diers2002} under the name of affine (or algebraic) set, but in a quite different context with no lattice-valued motivation or system notion. As a consequence our generalized topological spaces (resp., systems) are renamed affine spaces (resp., systems).


\subsection{Algebraic preliminaries}


\par
In this subsection, we recall the algebraic notions which will be used throughout the paper.

\begin{defn}
 \label{defn:1}
 Let $\Omega=(n_\lambda)_{\lambda\in\Lambda}$ be a family of cardinal numbers, which is indexed by a (possibly proper
 or empty) class $\Lambda$. An \emph{$\Omega$-algebra} is a pair $(A,(\omega^A_\lambda)_{\lambda\in\Lambda})$, which comprises a set $A$ and a family of maps $A^{n_\lambda} \arw{\omega^A_\lambda}A$ (\emph{$n_\lambda$-ary primitive operations} on $A$). An \emph{$\Omega$-homomorphism} $(A,(\omega^A_\lambda)_{\lambda\in \Lambda})\arw{\varphi}(B,(\omega^B_\lambda)_{\lambda\in\Lambda})$ is a map $A\arw{\varphi}B$, which makes the diagram
 $$
   \xymatrix{A^{n_\lambda} \ar[d]_-{\omega^A_\lambda} \ar[r]^-{\varphi^{n_\lambda}} & B^{n_\lambda}\ar[d]^-{\omega^B_\lambda}\\
             A \ar[r]_-{\varphi} & B\\}
 $$
 commute for every $\lambda\in\Lambda$. \algo\ is the construct of $\Omega$-algebras and $\Omega$-homomorphisms.
 \sqed
\end{defn}

\begin{defn}
 \label{defn:2}
 Let \mcal{M} (resp., \mcal{E}) be the class of $\Omega$-homomorphisms with injective (resp., surjective) underlying maps. A \emph{variety of $\Omega$-algebras} is a full subcategory of \algo, which is closed under the formation of products, \mcal{M}-subobjects (subalgebras), and \mcal{E}-quotients (homomorphic images). The objects (resp., morphisms) of a variety are called \emph{algebras} (resp., \emph{homomorphisms}).
 \sqed
\end{defn}

\par
In the following, we provide some examples of varieties, which are relevant to this paper.

\begin{exmp}
 \label{exmp:1}
 \hfill\par
 \begin{enumerate}[(1)]
  \item \semilatc{\bigvee} is the variety of \emph{$\bigvee$-semilattices}, i.e., partially ordered sets, which have arbitrary
        joins, and \semilatc{\bigwedge} is the variety of \emph{$\bigwedge$-semilattices}, i.e., partially ordered sets, which have arbitrary meets. $\semilatc{\bigvee}$ (resp. $\semilatc{\bigwedge}$) is precisely the category $\algop{\bigvee}$ (resp. $\algop{\bigwedge}$).
  \item \cat{Frm} is the variety of \emph{frames}, i.e., $\bigvee$-semilattices $A$, with singled out finite meets, and which
        additionally satisfy the distributivity condition $a\wedge(\bigvee S)=\bigvee_{s\in S}(a\wedge s)$ for every $a\in A$ and every $S\seq A$~\cite{Johnstone1982}. $\cat{Frm}$ is a full subcategory of $\algop{\bigvee,\wedge}$.
  \item \cat{CBAlg} is the variety of \emph{complete Boolean algebras}, i.e., complete lattices $A$ such that $a\wedge(b\vee
        c)=(a\wedge b)\vee(a\wedge c)$ for every $a$, $b$, $c\in A$, equipped with a unary operation $A\arw{(-)^{\ast}}A$ such that $a\vee a^{\ast}=\top_A$ and $a\wedge a^{\ast}=\bot_A$ for every $a\in A$, where $\top_A$ (resp., $\bot_A$) is the largest	(resp., smallest) element of $A$. \cat{CBAlg} is a full subcategory of \algop{\bigvee,\bigwedge,{}^{\ast}}.
  \item \cat{CSL} is the variety of \emph{closure semilattices}, i.e., $\bigwedge$-semilattices, with the singled out bottom
        element. $\cat{CSL}$ is a full subcategory of \algop{\bigwedge,\bot}.
        \sqed
 \end{enumerate}
\end{exmp}


\subsection{Affine spaces}


\par
In this subsection, we provide an extension of the notion of affine set of Y.~Diers~\cite{Diers1996,Diers1999,Diers2002}.

\begin{defn}
 \label{defn:3}
 Given a functor $\cat{X}\arw{T}\opm{\cat{B}}$, where \cat{B} is a variety of algebras, \topt{T} is the concrete category over \cat{X}, whose objects (\emph{$T$-affine spaces} or \emph{$T$-spaces}) are pairs $(X,\tau)$, where $X$ is an \cat{X}-object and $\tau$ is a \cat{B}-subalgebra of $TX$; and whose morphisms (\emph{$T$-affine morphisms} or \emph{$T$-morphisms}) $(X_1,\tau_1)\arw{f}(X_2,\tau_2)$ are \cat{X}-morphisms $X_1\arw{f}X_2$ with the property that $\opm{(Tf)}(\alpha)\in\tau_1$ for every $\alpha\in\tau_2$.
 \sqed
\end{defn}

\par
The following easy result will give rise to our main examples of $T$-spaces and $T$-morphisms.

\begin{prop}
 \label{prop:1}
 Given a variety \cat{B}, every subcategory \cat{S} of \loc{B} induces a functor $\set\times\cat{S}\arw{\mcal{P}_{\cat{S}}}\loc{B}$, $\mcal{P}_{\cat{S}}((X_1,B_1)\arw{(f,\varphi)}(X_2,B_2))=B_1^{X_1}\arw{\mcal{P}_{\cat{S}}(f,\varphi)}B_2^{X_2}$, where $\opm{(\mcal{P}_{\cat{S}}(f,\varphi))}(\alpha)=\opm{\varphi}\circ\alpha\circ f$.
\end{prop}

\par
The case $\cat{S}=\{B\arw{1_B}B\}$ provides a functor $\set\arw{\mcal{P}_{B}}\loc{B}$, $\mcal{P}_{B}(X_1\arw{f}X_2)=
B^{X_1}\arw{\mcal{P}_{B}f}B^{X_2}$, where $\opm{(\mcal{P}_{B}f)}(\alpha)=\alpha\circ f$. In particular, if $\cat{B}=\cat{CBAlg}$, and $\cat{S}=\{\msf{2}\arw{1_{\cat{2}}}\msf{2}\}$, then one obtains the well-known contravariant powerset functor $\set\arw{\mcal{P}}\loc{CBAlg}$, which is given on a map $X_1\arw{f}X_2$ by $\mcal{P}X_2\arw{\opm{(\mcal{P}f)}}\mcal{P}X_1$ with $\opm{(\mcal{P}f)}(S)=\{x\in X_1\,|\,f(x)\in S\}$.

\par
The following examples of the categories of the form \topt{T} will be relevant to this paper.

\begin{exmp}
 \label{exmp:2}
 \hfill\par
 \begin{enumerate}[(1)]
  \item If $\cat{B}=\cat{Frm}$, then $\topt{\mcal{P}_{\msf{2}}}$ is the category \cat{Top} of topological spaces.
  \item If $\cat{B}=\cat{CSL}$, then $\topt{\mcal{P}_{\msf{2}}}$ is the category \cat{Cls} of closure spaces~\cite{Aerts2002}.
  \item \topt{\mcal{P}_{B}} is the category $\cat{Af\,Set}(B)$ of affine sets of Y.~Diers.
  \item If $\cat{B}=\cat{Frm}$, then $\topt{\mcal{P}_{\cat{S}}}$ is the category \btop{\cat{S}} of variable-basis
        lattice-valued topological spaces of S.~E.~Rodabaugh~\cite{Rodabaugh1999a}.
        \sqed
 \end{enumerate}
\end{exmp}

\par
Given a map $X\arw{f}Y$ and subsets $R\seq X$, $S\seq Y$, $f^{\rightarrow}(R):=\{f(r)\,|\,r\in R\}$ and $f^{\leftarrow}(S):=\{x\,|\,f(x)\in S\}$~\cite{Rodabaugh1999c}; this notation will be employed through the rest of the paper.

\par
All of the categories of the form \topt{T} share the following convenient property.

\begin{thm}
 \label{thm:1}
 Given a functor $\cat{X}\arw{T}\opm{\cat{B}}$, the concrete category $(\topt{T},|-|)$ is topological (over \cat{X}).
\end{thm}
\begin{pf}
 Given a $|-|$-structured source $\mcal{S}=(X\arw{f_i}|(X_i,\tau_i)|)_{i\in I}$, its lift on $X$ w.r.t. \mcal{S} is given by the subalgebra of $TX$, which is generated by the union $\bigcup_{i\in I}(\opm{(Tf_i)})^{\rightarrow}(\tau_i)$;
 it is left to the reader to show that this lift is the unique initial lift. Alternatively, given a $|-|$-structured sink $\mcal{S}=(|(X_i,\tau_i)|\arw{f_i}X)_{i\in I}$, its lift on $X$ w.r.t. \mcal{S} is the intersection $\bigcap_{i\in I}(\opm{(Tf_i)})^{\leftarrow}(\tau_i)$; it is left to the reader to show that this lift is the unique final lift.
 \qed
\end{pf}

\par
As a consequence, one obtains the well-known result that all the categories of Example~\ref{exmp:2} are topological.


\subsection{Affine systems}


\par
Following the ideas of~\cite{Solovyov2012,Solovyov2013a}, this subsection introduces the concept of affine system as an analogue of topological systems of S.~Vickers~\cite{Vickers1989}.

\begin{defn}
 \label{defn:4}
 Given a functor $\cat{X}\arw{T}\loc{B}$, \topsyst{T} is the comma category $(T\downarrow 1_{\loc{B}})$, which is concrete
 over the product category $\cat{X}\times\loc{B}$, whose objects (\emph{$T$-affine systems} or \emph{$T$-systems}) are triples $(X,\kappa,B)$, which are made with \loc{B}-morphisms $TX\arw{\kappa}B$; and whose morphisms (\emph{$T$-affine morphisms} or \emph{$T$-morphisms}) $(X_1,\kappa_1,B_1)\arw{(f,\varphi)}(X_2,\kappa_2,B_2)$ are $\cat{X}\times\loc{B}$-morphisms $(X_1,B_1)\arw{(f,\varphi)}(X_2,B_2)$, which make the diagram
 $$
   \xymatrix{TX_1 \ar[d]_{\kappa_1} \ar[r]^{Tf} & TX_2 \ar[d]^{\kappa_2}\\
             B_1 \ar[r]_{\varphi} & B_2}
 $$
 commute.
 \sqed
\end{defn}

\begin{exmp}
 \label{exmp:3}
 \hfill\par
 \begin{enumerate}[(1)]
  \item If $\cat{B}=\cat{Frm}$, then \topsyst{\mcal{P}_{\msf{2}}} is the category \cat{TopSys} of topological systems of
        S.~Vickers.
  \item If $\cat{B}=\set$, then \topsyst{\mcal{P}_B} is the category $\cat{Chu}_B$ of Chu spaces over a set $B$ of
        P.-H.~Chu~\cite{Barr1979}.
  \item If $\cat{B}=\cat{Frm}$, then $\topsyst{\mcal{P}_{\cat{S}}}$ is the category \btopsys{\cat{S}} of variable-basis
        lattice-valued topological systems of J.~T.~Denniston, A.~Melton, and S.~E.~Rodabaugh~\cite{Denniston2009}.
        \sqed
 \end{enumerate}
\end{exmp}

\par
To provide another example of categories of the form \topsyst{T}, an additional notion is needed.

\begin{defn}
 \label{defn:5}
 A $T$-system $(X,\kappa,B)$ is called \emph{separated} provided that $TX\arw{\kappa}B$ is an epimorphism in \loc{B}. \spt{T} is the full subcategory of \topsyst{T} of separated $T$-systems.
 \sqed
\end{defn}

\par
We recall from, e.g.,~\cite[pp.~393~--~394]{Banaschewski1976} that monomorphisms in every variety \cat{B} are necessarily injective; given a monomorphism $B\arw{\varphi}\prm{B}$, the set $K=\{(b_1,b_2)\in B\times B\,|\,\varphi(b_1)=\varphi(b_2)\}$ is a subalgebra of $B\times B$ such that the respective projections \cell{K}{\pi_1}{\pi_2}{B} satisfy $\varphi\circ\pi_1=\varphi\circ\pi_2$, and therefore, $\pi_1=\pi_2$.

\begin{exmp}
 \label{exmp:4}
 If $\cat{B}=\cat{CSL}$, then \spt{\mcal{P}_{\msf{2}}} is the category \cat{SP} of state property systems of
 D.~Aerts~\cite{Aerts2002}.~\sqed
\end{exmp}


\subsection{Affine spatialization procedure}


\par
Following the results of~\cite{Solovyov2013a}, this subsection shows an affine analogue of the topological system spatialization procedure of S.~Vickers.

\begin{thm}
 \label{thm:2}
 \hfill\par
 \begin{enumerate}[(1)]
  \item There is a full embedding \incl{\topt{T}}{E}{\topsyst{T},} $E((X_1,\tau_1)\arw{f}(X_2,\tau_2))=
        (X_1,\opm{e}_{\tau_1},\tau_1)\arw{(f,\varphi)}(X_2,\opm{e}_{\tau_2},\tau_2)$, where $e_{\tau_i}$ is the inclusion $\tau_i\hookrightarrow TX_i$, and \opm{\varphi} is the restriction $\tau_2\arw{\opm{(Tf)}|_{\tau_2}^{\tau_1}}\tau_1$.
  \item $E$ has a right-adjoint-left-inverse $\topsyst{T}\arw{Spat}\topt{T}$, $Spat((X_1,\kappa_1,B_1)\arw{(f,\varphi)}
        (X_2,\kappa_2,B_2))=(X_1,(\opm{\kappa}_1)^{\rightarrow}(B_1))\arw{f}(X_2,(\opm{\kappa}_2)^{\rightarrow}(B_2))$.
  \item \topt{T} is isomorphic to a full (regular mono)-coreflective subcategory of \topsyst{T}.
 \end{enumerate}
\end{thm}
\begin{pf}
 To show that $Spat$ is a right adjoint to $E$, it is enough to verify that every system $(X,\kappa,B)$ has an $E$-co-universal arrow, i.e., a $T$-morphism $ESpat(X,\kappa,B)\arw{\varepsilon}(X,\kappa,B)$ such that for every $T$-morphism $E(\prm{X},\prm{\tau})\arw{(f,\varphi)}(X,\kappa,B)$, there exists a unique $T$-morphism $(\prm{X},\prm{\tau})\arw{g}Spat(X,\kappa,B)$ with $\varepsilon\circ Eg=(f,\varphi)$.

 \par
 There exists a $T$-morphism $(E Spat(X,\kappa,B)=(X,e^{op}_{(\opm{\kappa})^{\rightarrow}(B)},(\opm{\kappa})^{\rightarrow}(B)))
 \arw{\varepsilon=(1_X,\kappa)}(X,\kappa,B)$. Given a $T$-morphism $E(\prm{X},\prm{\tau})\arw{(f,\varphi)}(X,\kappa,B)$, it follows that $\opm{(Tf)}\circ\opm{\kappa}=e_{\prm{\tau}}\circ\opm{\varphi}$, which yields the desired $T$-morphism $(\prm{X},\prm{\tau})\arw{f}(Spat(X,\kappa,B)=(X,(\opm{\kappa})^{\rightarrow}(B)))$, whose uniqueness is clear.

 \par
 For the last claim, it is enough to show that given a $T$-system $(X,\kappa,B)$, the map $B\arw{\opm{\kappa}}(\opm{\kappa})^{\rightarrow}(B)$ is a regular epimorphism in \cat{B}. Define $C=\{(b_1,b_2)\in B\times B\,|\,\opm{\kappa}(b_1)=\opm{\kappa}(b_2)\}$ (the \emph{kernel} of $\opm{\kappa}$), and let $C\arw{\pi_i}B$ be given by $\pi_i(b_1,b_2)=b_i$ for $i\in\{1,2\}$. Then $(\opm{\kappa},(\opm{\kappa})^{\rightarrow}(B))$ is a coequalizer of $(\pi_1,\pi_2)$.
 \qed
\end{pf}

\par
The analogue of Theorem~\ref{thm:2} for the category \spt{T} is even better.

\begin{thm}
 \label{thm:3}
 The functors $E$ and $Spat$ restrict to \incl{\topt{T}}{\ovr{E}}{\spt{T}} and $\spt{T}\arw{\ovr{Spat}}\topt{T}$, respectively, providing an equivalence between the categories \topt{T} and \spt{T} such that $\ovr{Spat}\,\ovr{E}=1_{\topt{T}}$.
\end{thm}
\begin{pf}
 By Theorem~\ref{thm:2}, $\ovr{Spat}$ is a right-adjoint-left-inverse to $\ovr{E}$. To prove the theorem, it is enough to show that for every separated $T$-system $(X,\kappa,B)$, the $E$-co-universal arrow $ESpat(X,\kappa,B)\arw{\varepsilon=(1_X,\kappa)}(X,\kappa,B)$ from the proof of Theorem~\ref{thm:2} is an isomorphism. The claim follows from the definition of $\varepsilon$, since $B\arw{\opm{\kappa}}(\opm{\kappa})^{\rightarrow}(B)$ is always surjective, and it is injective by the property of separated $T$-systems.
 \qed
\end{pf}

\begin{cor}
 \label{cor:1}
 The category \topt{T} is the amnestic modification of the category \spt{T}.
\end{cor}
\begin{pf}
 Follows from Theorem~\ref{thm:3} and the definition of the amnestic modification of~\cite[Remark~5.34]{Adamek2009}.
 \qed
\end{pf}

\par
The following (well-known) results are direct consequences of Theorems~\ref{thm:2} and \ref{thm:3}, respectively.

\begin{rem}
 \label{rem:1}
 \hfill\par
 \begin{enumerate}[(1)]
  \item \cat{Top} is isomorphic to a full (regular mono)-coreflective subcategory of the category \cat{TopSys}, which provides
        the system spatialization procedure of S.~Vickers. More generally, \btop{\cat{Loc}} (\cat{Loc} is the dual of \cat{Frm}) is isomorphic to a full (regular mono)-coreflective subcategory of the category \btopsys{\cat{Loc}}.
  \item The categories \cat{Cls} and \cat{SP} are equivalent~\cite{Aerts1999a,Aerts2002}.
        \sqed
 \end{enumerate}
\end{rem}

\par
Moreover, from Corollary~\ref{cor:1}, one gets~\cite[Theorem~4]{Aerts2002}, which states that the category \cat{Cls} is the amnestic modification of the category \cat{SP}.


\section{Localification procedure for affine systems}


\par
This section provides a localification procedure for affine systems, motivated by the above-mentioned localification procedure for topological systems of S.~Vickers.

\begin{prop}
 \label{prop:2}
 There is a functor $\topsyst{T}\arw{Loc}\loc{B}$, $Loc((X_1,\kappa_1,B_1)\!\arw{(f,\varphi)}\!(X_2,\kappa_2,B_2))=
 B_1\arw{\varphi}B_2$.
\end{prop}

\par
Unlike the affine spatialization procedure, in which the functor in the opposite direction always exists, the localification procedure is more demanding.

\begin{thm}
 \label{thm:4}
 Given a functor $\cat{X}\arw{T}\loc{B}$, the following are equivalent.
 \begin{enumerate}[(1)]
  \item There exists an adjoint situation $(\eta,\varepsilon):T\dashv Pt:\loc{B}\arw{}\cat{X}$.
  \item There exists a full embedding \incl{\loc{B}}{E}{\topsyst{T}} such that $Loc$ is a left-adjoint-left-inverse to $E$.
        \loc{B} is then isomorphic to a full reflective subcategory of \topsyst{T}.
 \end{enumerate}
\end{thm}
\begin{pf}
 \hfill\par
 \emph{Ad} $(1)\Rightarrow(2)$. Define a functor $\loc{B}\arw{E}\topsyst{T}$ by $E(B_1\arw{\varphi}B_2)=(Pt B_1,\varepsilon_{B_1},B_1)\arw{(Pt\varphi,\varphi)}(Pt B_2,\varepsilon_{B_2},B_2)$. Correctness of $E$ on morphisms follows from commutativity of the diagram
 $$
   \xymatrix{TPt B_1\ar[rr]^{TPt\varphi} \ar[d]_-{\varepsilon_{B_1}} && TPt B_2 \ar[d]^-{\varepsilon_{B_2}}\\
             B_1 \ar[rr]_-{\varphi} && B_2.}
 $$
 Moreover, $E$ is clearly an embedding. To verify that the functor $E$ is full, note that given a $T$-morphism $(Pt B_1,\varepsilon_{B_1},B_1)\arw{(f,\varphi)}(Pt B_2,\varepsilon_{B_2},B_2)$, commutativity of the diagrams
 $$
   \xymatrix{TPt B_1\ar[rr]^-{\varepsilon_{B_1}} \ar@<-1ex>[d]_-{TPt\varphi}	
             \ar@<1ex>[d]^-{Tf} && B_1 \ar[d]^-{\varphi}\\
             TPt B_2 \ar[rr]_-{\varepsilon_{B_2}} && B_2}
 $$
 implies that $\varepsilon_{B_2}\circ TPt\varphi=\varepsilon_{B_2}\circ Tf$, and therefore, $Pt\varphi=f$. Given a $T$-system $(X,\kappa,B)$, straightforward calculations show that $(X,\kappa,B)\arw{(f:=Pt\kappa\circ\eta_X,1_B)}((PtB,\varepsilon_B,B)=ELoc(X,\kappa,B))$ provides an $E$-universal arrow for $(X,\kappa,B)$. It is also easy to see that $Loc E=1_{\loc{B}}$.

 \par
 \emph{Ad} $(2)\Rightarrow(1)$. Given an adjunction $Loc\dashv E:\loc{B}\arw{}\topsyst{T}$, $\cat{X}\arw{T}\loc{\cat{B}}$ is the composition of the left adjoint functors $\cat{X}\arw{}\topt{T}$ (the \emph{indiscrete functor} of, e.g.,~\cite[Proposition~21.12\,(2)]{Adamek2009}, which exists by Theorem~\ref{thm:1}), \incl{\topt{T}}{E}{\topsyst{T}} (the embedding of Theorem~\ref{thm:2}), and $\topsyst{T}\!\arw{Loc}\!\opm{\cat{B}}$.
 \qed
\end{pf}

\par
The following provides an example of the functor $Pt$ of Theorem~\ref{thm:4}\,(1).

\begin{thm}
 \label{thm:5}
 Every functor $\set\arw{\mcal{P}_B}\loc{B}$ has a right adjoint.
\end{thm}
\begin{pf}
 Recall that the right adjoint functor $\loc{B}\arw{Pt_B}\set$ is given by $Pt_B(B_1\arw{\varphi}B_2)=
 \cat{B}(B_1,B)\arw{Pt_B\varphi}\cat{B}(B_2,B)$, where $(Pt_B\varphi)(p)=p\circ\opm{\varphi}$. Given a \cat{B}-algebra $A$, the map $A\arw{\opm{\varepsilon}}(\mcal{P}_BPt_BA=B^{\cat{B}(A,B)})$, defined by $(\opm{\varepsilon}(a))(p)=p(a)$, provides a $\mcal{P}_B$-co-universal arrow for $A$.
 \qed
\end{pf}

\par
As a consequence of Theorems~\ref{thm:4}, \ref{thm:5}, one gets the following well-known results.

\begin{rem}
 \label{rem:2}
 \hfill\par
 \begin{enumerate}[(1)]
  \item \cat{Loc} (the dual of \cat{Frm}) is isomorphic to a full reflective subcategory of \cat{TopSys}, which provides the
        system localification procedure of S.~Vickers.
  \item \loc{B} is isomorphic to a full reflective subcategory of \topsyst{\mcal{P}_B}.
        \sqed
 \end{enumerate}
\end{rem}

\par
In particular, the case of the category \cat{TopSys} shows that in Theorem~\ref{thm:4}\,(2), the category \loc{B}, even if a
reflective subcategory of \topsyst{T}, can be both non-mono-reflexive and non-epi-reflexive.

\par
In~\cite{Rodabaugh1999a}, S.~E.~Rodabaugh considered functors of the form $\set\times\cat{S}\arw{\mcal{P}_{\cat{S}}}\cat{Loc}$ and their respective categories of affine spaces, using, however, a different terminology (recall Example~\ref{exmp:2}\,(4)). The next result shows that Theorem~\ref{thm:5} (in general) can not be extended from the subcategory $\{B\arw{1_B}B\}$ to the whole \loc{\cat{B}}.

\begin{prop}
 \label{prop:3}
 Consider a functor $\set\times\loc{B}\arw{T:=\mcal{P}_{\loc{B}}}\loc{B}$, and suppose that there exists a \cat{B}-algebra $B$, whose underlying set has finite cardinality $n > 1$. Then $T$ has no right adjoint.
\end{prop}
\begin{pf}
 If $T$ has a right adjoint, then $T$ preserves coproducts. Given a singleton set \msf{1}, $T((\msf{1},B)\coprod(\msf{1},B))=T((\msf{1}\biguplus\msf{1},B\times B))=(B\times B)^{(\msf{1}\biguplus\msf{1})}$ and $T(\msf{1},B)\times T(\msf{1},B)=B^{\msf{1}}\times B^{\msf{1}}$. Since $T((\msf{1},B)\coprod(\msf{1},B))\cong T(\msf{1},B)\times T(\msf{1},B)$, $n^4=Card((B\times B)^{(\msf{1}\biguplus\msf{1})})=Card(B^{\msf{1}}\times B^{\msf{1}})=n^2$, which is a contradiction.
 \qed
\end{pf}

\par
For instance, Proposition~\ref{prop:3} implies that the functor $\set\times\cat{Loc}\arw{\mcal{P}_{\cat{Loc}}}\cat{Loc}$ has no right adjoint, i.e., Theorem~\ref{thm:4}\,(2) is not applicable to the category \btopsys{\cat{Loc}} of Example~\ref{exmp:3}\,(3).


\subsection{Affine theories}


\par
The results of this subsection stem from~\cite{Solovyovh}, and are motivated by the categorical approach to universal algebra, based on the concept of \emph{algebraic theory} of~F.~W.~Lawvere~\cite{Lawvere1963}. We would like to develop the tools, which would allow one to study the properties of a category of the form \topsyst{T} (or \topt{T}) through the properties of its related functor $\cat{X}\arw{T}\opm{\cat{B}}$.

\begin{defn}
 \label{defn:6}
 An \emph{affine theory} is a functor $\cat{X}\arw{T}\opm{\cat{B}}$, where \cat{B} is a variety of algebras.
 \sqed
\end{defn}

\par
Throughout this paper, there is no distinction made between categories and quasicategories~\cite{Adamek2009}. The latter are defined similarly to the former except that their objects do not necessarily form a class and their hom-families are not necessarily sets.

\begin{defn}
 \label{defn:7}
 \cat{AfTh} is the category, whose objects are affine theories $\cat{X}\arw{T}\loc{B}$, and whose morphisms $T_1\arw{(F,\Phi,\eta)}T_2$ (also shortened to $\eta$) comprise two functors $\cat{X}_1\arw{F}\cat{X}_2$, $\cat{B}_1\arw{\Phi}\cat{B}_2$ and a natural transformation $T_2F\arw{\eta}\loc{$\Phi$}T_1$, or, more specifically,
 $$
   \xymatrix{\cat{X}_1 \ar[d]_{T_1}
             \ar[rr]^{F} & & \cat{X}_2 \ar@{=>}"2,1"_-{\eta}
             \ar[d]^{T_2}\\
             \loc{B}_1 \ar[rr]_{\opm{\Phi}} & & \loc{B}_2.}
 $$
 Given two affine theories $T_1\arw{\eta_1}T_2$ and $T_2\arw{\eta_2}T_3$, their composition is defined by $T_3F_2F_1\arw{\eta_2\odot\eta_1}\opm{\Phi}_2\opm{\Phi}_1T_1=T_3F_2F_1\arw{{\eta_2}_{F_1}}\opm{\Phi}_2T_2F_1\arw{\opm{\Phi}_2\eta_1}
 \opm{\Phi}_2\opm{\Phi}_1T_1$. The identity on an affine theory $T$ is provided by the  identity natural transformation $T\arw{1_T}T$.
 \sqed
\end{defn}

\par
Note that the functor $\opm{\Phi}$ in Definition~\ref{defn:7} goes in the same direction as $\Phi$, i.e., does not switch its domain and codomain. The next lemma shows that the construction of Definition~\ref{defn:7} is indeed a category.

\begin{lem}
 \label{lem:1}
 The composition law is associative and the identities are as given.
\end{lem}
\begin{pf}
 Consider three affine theory morphisms $T_1\arw{\eta_1}T_2$, $T_2\arw{\eta_2}T_3$ and $T_3\arw{\eta_3}T_4$, or, more specifically, the following diagram:
 $$
   \xymatrix{\cat{X}_1 \ar[d]_{T_1}
             \ar[rr]^{F_1} & & \cat{X}_2 \ar@{=>}"2,1"_-{\eta_1}
             \ar[d]^{T_2} \ar[rr]^{F_2} & & \cat{X}_3 \ar[d]^{T_3} \ar@{=>}"2,3"_-{\eta_2} \ar[rr]^{F_3} & &  \cat{X}_4 \ar[d]^{T_4} \ar@{=>}"2,5"_-{\eta_3}\\
             \loc{B}_1 \ar[rr]_{\opm{\Phi}_1} & & \loc{B}_2 \ar[rr]_{\opm{\Phi}_2} & & \loc{B}_3 \ar[rr]_{\opm{\Phi}_3}  && \loc{B}_4.}
 $$
 Straightforward computations show that
 \begin{align*}
  T_4F_3F_2F_1\arw{(\eta_3\odot\eta_2)\odot\eta_1}\opm{\Phi}_3
  \opm{\Phi}_2\opm{\Phi}_1T_1=\\
  T_4F_3F_2F_1\arw{(\eta_3\odot\eta_2)_{F_1}}\opm{\Phi}_3
  \opm{\Phi}_2T_2F_1\arw{\opm{\Phi}_3\opm{\Phi}_2\eta_1}\opm{\Phi}_3
  \opm{\Phi}_2\opm{\Phi}_1T_1=\\
  T_4F_3F_2F_1\arw{{\eta_3}_{F_2F_1}}\opm{\Phi}_3T_3F_2F_1
  \arw{\opm{\Phi}_3{\eta_2}_{F_1}}\opm{\Phi}_3\opm{\Phi}_2T_2F_1
  \arw{\opm{\Phi}_3\opm{\Phi}_2\eta_1}\opm{\Phi}_3
  \opm{\Phi}_2\opm{\Phi}_1T_1=\\
  T_4F_3F_2F_1\arw{{\eta_3}_{F_2F_1}}
  \opm{\Phi}_3(T_3F_2F_1\arw{{\eta_2}_{F_1}}
  \opm{\Phi}_2T_2F_1\arw{\opm{\Phi}_2\eta_1}\opm{\Phi}_2\opm{\Phi}_1T_1)=\\
  T_4F_3F_2F_1\arw{{\eta_3}_{F_2F_1}}\opm{\Phi}_3T_3F_2F_1
  \arw{\opm{\Phi}_3(\eta_2\odot\eta_1)}
  \opm{\Phi}_3\opm{\Phi}_2\opm{\Phi}_1T_1=\\
  T_4F_3F_2F_1\arw{\eta_3\odot(\eta_2\odot\eta_1)}\opm{\Phi}_3
  \opm{\Phi}_2\opm{\Phi}_1T_1.\hspace{4mm}
 \end{align*}
 The last statement of the lemma is clear.
 \qed
\end{pf}

\par
The composition law of the category \cat{AfTh} resembles the \emph{star product} of \cite[Definition~13.10]{Herrlich2007}, but does not coincide with it.

\begin{defn}
 \label{defn:8}
 \cat{AfStm} is the category, whose objects are categories of the form \topsyst{T} and whose morphisms are functors between them.
 \sqed
\end{defn}

\par
The following theorem shows that affine theory morphisms may be translated into functors between their respective categories of affine systems.

\begin{thm}
 \label{thm:6}
 There is a functor $\cat{AfTh}\!\arw{AfSys}\!\cat{AfStm}$, $AfSys(T_1\!\arw{\eta}\!T_2)\!=\!\topsyst{T_1}\!\arw{AfSys\eta}\!
 \topsyst{T_2}$, where $AfSys\eta((X,\kappa,B)\arw{(f,\varphi)}(\prm{X},\prm{\kappa},\prm{B}))=
 (FX,\opm{\Phi}\kappa\circ\eta_X,\opm{\Phi}B)\arw{(Ff,\opm{\Phi}\varphi)}
 (F\prm{X},\opm{\Phi}\prm{\kappa}\circ\eta_{\prm{X}},\opm{\Phi}\prm{B})$.
\end{thm}
\begin{pf}
 We begin with the verification that the definition of $AfSys\eta$ is correct. Given an \topsyst{T_1}-morphism $(X,\kappa,B)\arw{(f,\varphi)}(\prm{X},\prm{\kappa},\prm{B})$, the diagram
 $$
   \xymatrix{T_1X\ar[d]_{\kappa}
             \ar[r]^{T_1f} & T_1\prm{X}
             \ar[d]^{\prm{\kappa}}\\
             B \ar[r]_{\varphi} & \prm{B}}
 $$
 commutes and that provides the commutativity of the next one:
 $$
   \xymatrix{T_2FX\ar[d]_{\eta_X}
             \ar[rr]^{T_2Ff} & & T_2F\prm{X}
             \ar[d]^{\eta_{\prm{X}}}\\
             \opm{\Phi}T_1X \ar[d]_{\opm{\Phi}\kappa} \ar[rr]^{\opm{\Phi}T_1f} & & \opm{\Phi}T_1\prm{X}
             \ar[d]^{\opm{\Phi}\prm{\kappa}}\\
             \opm{\Phi}B \ar[rr]_{\opm{\Phi}\varphi} && \opm{\Phi}
             \prm{B}.}
 $$
 It follows that $(FX,\opm{\Phi}\kappa\circ\eta_X,\opm{\Phi}B)\arw{(Ff,\opm{\Phi}\varphi)}
 (F\prm{X},\opm{\Phi}\prm{\kappa}\circ\eta_{\prm{X}},\opm{\Phi}\prm{B})$ is an \topsyst{T_2}-morphism.

 \par
 To show that $AfSys$ also preserves the composition, we notice that given two affine theory morphisms $T_1\arw{(F_1,\Phi_1,\eta_1)}T_2$ and $T_2\arw{(F_2,\Phi_2,\eta_2)}T_3$, it follows that
 \begin{align*}
  AfSys((F_2,\Phi_2,\eta_2)\circ(F_1,\Phi_1,\eta_1))((X,\kappa,B)\arw{(f,\varphi)}(\prm{X},\prm{\kappa},\prm{B}))=\\
  AfSys(F_2F_1,\Phi_2\Phi_1,\eta_2\odot\eta_1)((X,\kappa,B)\arw{(f,\varphi)}(\prm{X},\prm{\kappa},\prm{B})))=\\
  (F_2F_1X,\opm{\Phi}_2\opm{\Phi}_1\kappa\circ(\eta_2\odot\eta_1)_X,\opm{\Phi}_2\opm{\Phi}_1B,)\arw{(F_2F_1f,
  \opm{\Phi}_2\opm{\Phi}_1\varphi)}\hspace{4mm}\\
  (F_2F_1\prm{X},\opm{\Phi}_2\opm{\Phi}_1\prm{\kappa}\circ(\eta_2\odot\eta_1)_{\prm{X}},
  \opm{\Phi}_2\opm{\Phi}_1\prm{B})\overset{(\dagger)}{=}\\
  (F_2F_1X,\opm{\Phi}_2(\opm{\Phi}_1\kappa\circ{\eta_1}_X)\circ{\eta_2}_{F_1X},\opm{\Phi}_2\opm{\Phi}_1B)
  \arw{(F_2F_1f,\opm{\Phi}_2\opm{\Phi}_1\varphi)}\hspace{4mm}\\
  (F_2F_1\prm{X},\opm{\Phi}_2(\opm{\Phi}_1
  \prm{\kappa}\circ{\eta_1}_{\prm{X}})\circ{\eta_2}_{F_1\prm{X}},
  \opm{\Phi}_2\opm{\Phi}_1\prm{B})=\\
  AfSys(F_2,\Phi_2,\eta_2)((F_1X,\opm{\Phi}_1\kappa\circ{\eta_1}_X,\opm{\Phi}_1B)\arw{(F_1f,\opm{\Phi}_1\varphi)}\hspace{4mm}\\
  (F_1\prm{X},\opm{\Phi}_1\prm{\kappa}\circ{\eta_1}_{\prm{X}},\opm{\Phi}_1\prm{B}))=\\
  AfSys(F_2,\Phi_2,\eta_2)\circ AfSys(F_1,\Phi_1,\eta_1)((X,\kappa,B)
  \arw{(f,\varphi)}(\prm{X},\prm{\kappa},\prm{B})),\hspace{4mm}
 \end{align*}
 where $(\dagger)$ relies on the fact that $\opm{\Phi}_2\opm{\Phi}_1\kappa\circ (\eta_2\odot\eta_1)_X=\opm{\Phi}_2\opm{\Phi}_1\kappa\circ\opm{\Phi}_2{\eta_1}_X\circ{\eta_2}_{F_1X}=\opm{\Phi}_2(\opm{\Phi}_1\kappa
 \circ{\eta_1}_X)\circ{\eta_2}_{F_1X}$.

 \par
 Preservation of the identities is straightforward.
 \qed
\end{pf}

\par
It is an interesting and challenging open question, whether the functor of Theorem~\ref{thm:6} has a left or right adjoint.

\par
For the sake of completeness, it is noted that by analogy with the category \cat{AfStm}, one can define the category, which comprises the categories of the form \topt{T}. The respective functor from \cat{AfTh} in this case though requires more effort and will be not considered here; this topic is studied in full detail in~\cite{Solovyovh}.


\section{Affine systems as a framework for elementary institutions}


\par
In this section is proposed a new framework for doing a certain type of institutions, namely, elementary institutions.  This new framework is based in affine systems.


\subsection{Institutions and their morphisms}


\par
The definitions of institution and its morphism are recalled from, e.g.,~\cite{Goguen1992} (elementary institutions or also pre-institutions come from~\cite{Salibra1993,Sernadas1995}). From now on, \cat{Cat} denotes the category of categories (recall there is no distinction made between categories and quasicategories).

\begin{defn}
 \label{defn:9}
 An \emph{institution} \mbb{I} consists of:
 \begin{itemize}
  \item a category \cat{Sign} of \emph{signatures} with $\Sigma$ denoting an arbitrary object,
  \item a functor $\cat{Sign}\arw{Mod}\opm{\cat{Cat}}$ giving \emph{$\Sigma$-models} and \emph{$\Sigma$-morphisms},
  \item a functor $\cat{Sign}\arw{Sen}\cat{Set}$ giving \emph{$\Sigma$-sentences},
  \item a \emph{satisfaction} relation $\models_{\Sigma}\seq\obc{Mod\Sigma}\times Sen\Sigma$ for every
        $\Sigma\in\obc{\cat{Sign}}$
 \end{itemize}
 such that for every \cat{Sgn}-morphism $\Sigma\arw{\phi}\prm{\Sigma}$, the \emph{Satisfaction Condition}
 $$
   \prm{m}\models_{\prm{\Sigma}}Sen\phi(s)\text{ iff }Mod\phi(\prm{m})\models_{\Sigma}s
 $$
 holds for every $\prm{m}\in\obc{Mod\prm{\Sigma}}$ and every $s\in Sen\Sigma$. An institution is called \emph{elementary}~\cite[Definition~3.1]{Sernadas1994},~\cite{Sernadas1995} (also~\emph{pre-institution}~\cite[Definition~1]{Salibra1993}) provided that the category \cat{Cat} (in the definition of the functor $Mod$) is replaced with the category \cat{Set}.
 \sqed
\end{defn}

\par
Elementary institutions are not much popular among the researchers at the moment. There do exist though several applications of this concept, which could be found in, e.g.,~\cite{Ehrich2000,Salibra1993,Sernadas1994,Sernadas1995}.

\par
We notice that~\cite{Goguen1986} defines elementary institutions under the name of \emph{simplest institution}. There also exists a more sophisticated concept of elementary institution of~\cite[p.~68]{Diaconescu2008}, which is not used in this paper.

\par
The above satisfaction relation $\models_{\Sigma}\seq\obc{Mod\Sigma}\times Sen\Sigma$ can be rewritten as a map
$\obc{Mod\Sigma}\arw{\kappa}\mcal{P}(Sen\Sigma)$, in which $\kappa(m)=\{e\in Sen\Sigma\,|\,m\models_{\Sigma}e\}$. Moreover, if one defines \cat{2} as the category $\bot\arw{\iota}\top$ (two objects and one non-identity morphism), every functor $Mod\Sigma\arw{K}\cat{2}^{Sen\Sigma}$ will provide a satisfaction relation $\models_{\Sigma}\seq\obc{Mod\Sigma}\times Sen\Sigma$ with $m\models_\Sigma e$ iff $Km(e)=\top$.

\begin{defn}
 \label{defn:10}
 Let \mbb{I} and $\prm{\mbb{I}}$ be institutions. Then an \emph{institution morphism} $\mbb{I}\arw{(\Phi,\alpha,\beta)}\prm{\mbb{I}}$ consists of:
 \begin{itemize}
  \item a functor $\cat{Sign}\arw{\Phi}\prm{\cat{Sign}}$,
  \item a natural transformation $\prm{Sen}\Phi\arw{\alpha}Sen$, or, more specifically,
        $$
          \xymatrix{\cat{Sign} \ar[d]_-{Sen} \ar[rr]^-{\Phi} & & \prm{\cat{Sign}}
                    \ar@{=>}"2,1"_-{\alpha} \ar[d]^{\prm{Sen}}\\
                    \cat{Set} \ar@{=}[rr] & & \cat{Set}}
        $$
  \item a natural transformation $Mod\arw{\beta}\prm{Mod}\Phi$ or, more specifically,
        $$
          \xymatrix{\cat{Sign} \ar[d]_-{Mod} \ar[rr]^-{\Phi} & & \prm{\cat{Sign}}
                    \ar[d]^{\prm{Mod}}\\
                    \opm{\cat{Cat}} \ar@{=>}"1,3"^-{\beta} \ar@{=}[rr] & & \opm{\cat{Cat}}}
        $$
 \end{itemize}
 such that the following \emph{satisfaction condition} holds
 $$
   m\models_{\Sigma}\alpha_{\Sigma}(\prm{s})\text{ iff }\beta_{\Sigma}(m)\prm{\models}_{\Phi\Sigma}\prm{s}
 $$
 for every $\Sigma$-model $m$ from \mbb{I} and every $\Phi\Sigma$-sentence $\prm{s}$ from $\prm{\mbb{I}}$.
 \sqed
\end{defn}

\par
Definitions~\ref{defn:9}, \ref{defn:10} provide the category \cat{Inst} (resp., \cat{ElInst}) of (resp., elementary)
institutions and their morphisms. Examples of institutions include, e.g., first-order logic (with equality), Horn clause logic (with equality), equational logic, order-sorted equational logic, continuous equational logic~\cite{Goguen1986} (we again recommend the reader to have a look at the more comprehensive list of examples of~\cite[Subsection~3.2]{Diaconescu2008}).

\par
Following~\cite{Goguen1992}, the direction of natural transformations $\alpha$, $\beta$ in Definition~\ref{defn:10} is (in a certain sense) arbitrary. More precisely, having natural transformations $Sen\arw{\alpha}\prm{Sen}\Phi$ and $\prm{Mod}\Phi\arw{\beta}Mod$ in hand, one gets the satisfaction condition
$$
  \beta_{\Sigma}(\prm{m})\models_{\Sigma}s\text{ iff }\prm{m}\prm{\models}_{\Phi\Sigma}\alpha_{\Sigma}(s)
$$
for every $\Sigma$-sentence $s$ from \mbb{I} and every $\Phi\Sigma$-model $\prm{m}$ from $\prm{\mbb{I}}$. Such a structure is called \emph{institution comorphism}~\cite{Goguen2002}. We notice that comorphisms are more used than morphisms in institution theory. Moreover, the duality between morphisms and comorphisms is explained in~\cite{Diaconescu2008}.


\subsection{Topological institutions and their morphisms}


\par
Following~\cite{Sernadas1995}, we recall the notion of topological institution, which is based in the category \cat{TopSys} of topological systems of S.~Vickers (recall Example~\ref{exmp:3}\,(1)).

\begin{defn}
 \label{defn:11}
 A \emph{topological institution} consists of:
 \begin{itemize}
  \item a category \cat{Sign} (of signatures),
  \item a functor $\cat{Sign}\arw{\mcal{T}}\opm{\cat{TopSys}}$.
        \sqed
 \end{itemize}
\end{defn}

\begin{defn}
 \label{defn:12}
 A \emph{topological institution morphism} $(\cat{Sign},\mcal{T})\arw{(\Phi,\alpha)}(\prm{\cat{Sign}},\prm{\mcal{T}})$ consists of:
 \begin{itemize}
  \item a functor $\cat{Sign}\arw{\Phi}\prm{\cat{Sign}}$,
  \item a natural transformation $\mcal{T}\arw{\alpha}\prm{\mcal{T}}\Phi$.
        \sqed
 \end{itemize}
\end{defn}

\par
Definitions~\ref{defn:11}, \ref{defn:12} provide the category \cat{TpInst} of topological institutions and their morphisms. Moreover, in~\cite{Sernadas1994,Sernadas1995} an adjoint situation between the categories \cat{ElInst} and \cat{TpInst} is constructed. For convenience of the reader, we provide a brief description of some of its details.

\par
We recall first from~\cite{Vickers1989} that a \emph{topological system} is a triple $(X,A,\models)$ (denoted $D=(\pt D,\omg D,\models)$), where $X$ is a set, $A$ is a locale, and $\models\seq X\times A$ is a binary relation such that

\begin{enumerate}[(1)]
 \item if $S$ is a finite subset of $A$, then for every $x\in X$, $x\models\bigwedge A$ iff $x\models a$ for every $a\in A$;
 \item if $S$ is any subset of $A$, then for every $x\in X$, $x\models\bigvee A$ iff $x\models a$ for some $a\in A$.
\end{enumerate}

\noindent
A \emph{topological system morphism} $D_1\arw{f=(\pt f,\,\omg f)}D_2$ consists of a map $\pt D_1\arw{\pt f}\pt D_2$ and a homomorphism of locales $\omg D_1\arw{\omg f}\omg D_2$ such that for every $x_1\in \pt D_1$, $a_2\in \omg D_2$, $x_1\models_1\opm{(\omg f)}(a_2)$ iff $\pt f(x_1)\models_2 a_2$.

\par
We continue with some necessary institutional preliminaries from~\cite{Sernadas1994,Sernadas1995}. Let $(\cat{Sign},\,Mod,\,Sen,\,\models)$ be an elementary institution, and let $\Sigma\in\obc{\cat{Sign}}$ be an abstract signature. If $\Phi\seq Sen\Sigma$ and $\psi\in Sen\Sigma$, then $\Phi$ \emph{semantically entails} $\psi$ (denoted $\Phi\vdash_{\Sigma}\psi$) provided that for every $x\in Mod\Sigma$, $x\models_{\Sigma}\varphi$ for every $\varphi\in\Phi$, implies $x\models_{\Sigma}\psi$. One denotes $\Phi^{\vdash_{\Sigma}}=\{\psi\in Sen\Sigma\mid\Phi\vdash_{\Sigma}\psi\}$. $\Phi$ is then called a \emph{theory} provided that $\Phi^{\vdash_{\Sigma}}=\Phi$. The class of all theories over $\Sigma$ is denoted $|Th_{\Sigma}|$. On the class $|Th_{\Sigma}|$, one introduces a partial order by $\Phi_1\leqs_{\Sigma}\Phi_2$ iff $\Phi_1\vdash_{\Sigma}\Phi_2$ (in the sense of the above semantical entailment relation $\vdash_{\Sigma}$) and obtains a complete lattice $Th^{op}_{\Sigma}=(|Th_{\Sigma}|,\leqs_{\Sigma})$, in which, given $S\seq|Th_{\Sigma}|$, $\bigvee S=\bigcap S$ and $\bigwedge S=(\bigcup S)^{\vdash_{\Sigma}}$. We get thus a triple $(Mod\Sigma,\,Th^{op}_{\Sigma},\,\Vdash_{\Sigma})$, in which $x\Vdash_{\Sigma}\Phi$ means $x\models_{\Sigma}\varphi$ for every $\varphi\in\Phi$. One can use now the coverage technique of~\cite{Johnstone1982} to complete the $\bigwedge$-semilattice of theories, while preserving only the joins, which are respected by the satisfaction relation. Afterwards one extends the satisfaction relation to the new theories. In such a way, one obtains a functor $\cat{ElInst}\arw{Top}\cat{TpInst}$ (we omit its definition on morphisms).

\par
The functor in the opposite direction $\cat{TpInst}\arw{Geo}\cat{ElInst}$ is easier. Given a topological institution $\mcal{T}$, one obtains an elementary institution $(\cat{Sign},\,Mod,\,Sen,\,\models)$, in which $Mod\Sigma=\pt \mcal{T}\Sigma$, $Sen\Sigma=|\omg\mcal{T}\Sigma|$ (where $|-|$ stands for the underlying set), and $\models_{\Sigma}\,=\,\models_{\mcal{T}\Sigma}$. One can show that $Geo$ is a left adjoint to $Top$ (\cite[Theorem~4.44]{Sernadas1994},\cite[Theorem~34]{Sernadas1995}). It is not clear to us at the moment, whether one could provide a generalization of the just mentioned machinery for the category \cat{Inst} of institutions.


\subsection{Affine institutions and their morphisms}


\par
Using the ideas of the previous subsection, the concepts of affine institution and affine institution morphism are now introduced.

\begin{defn}
 \label{defn:13}
 An \emph{affine institution} consists of:
 \begin{itemize}
  \item a category \cat{S} (of abstract signatures),
  \item an affine theory $T$,
  \item a functor $\cat{S}\arw{I}\topsyst{T}$.
        \sqed
 \end{itemize}
\end{defn}

\par
Following the remark of the referee, we notice that unlike the case of topological institutions, affine institutions have their category of signatures denoted \cat{S}. The only reason for the change is our wish to, first, shorten the notations related to affine institutions, and, second, underline the \emph{abstract} nature of signatures.

\begin{defn}
 \label{defn:14}
 An \emph{affine institution morphism} $(\cat{S}_1,T_1,I_1)\arw{(\Phi,\eta,\alpha)}(\cat{S}_2,T_2,I_2)$ consists of:
 \begin{itemize}
  \item a functor $\cat{S}_1\arw{\Phi}\cat{S}_2$,
  \item an affine theory morphism $T_1\arw{\eta}T_2$,
  \item a natural transformation $AfSys\eta I_1\arw{\alpha}I_2\Phi$ (recall the notation of Theorem~\ref{thm:6}), or, more
        specifically,
        $$
          \xymatrix{\cat{S}_1 \ar[d]_-{I_1} \ar[rr]^-{\Phi} & & \cat{S}_2 \ar[d]^{I_2}\\
                    \topsyst{T_1} \ar@{=>}"1,3"^-{\alpha} \ar[rr]_-{AfSys\eta} & & \topsyst{T_2}.}
        $$
        \sqed
 \end{itemize}
\end{defn}

\par
Definitions~\ref{defn:13}, \ref{defn:14} provide the category \cat{AfInst} of affine institutions and their morphisms. For example, given affine institution morphisms
$$
  \xymatrix{\cat{S}_1 \ar[d]_-{I_1} \ar[rr]^-{\Phi} & & \cat{S}_2 \ar[d]^{I_2} \ar[rr]^-{\prm{\Phi}} & & \cat{S}_3  \ar[d]^-{I_3}\\
            \topsyst{T_1} \ar@{=>}"1,3"^-{\alpha} \ar[rr]_-{AfSys\eta} & & \topsyst{T_2} \ar@{=>}"1,5"^-{\prm{\alpha}}\ar[rr]_-{AfSys\prm{\eta}} && \topsyst{T_3},}
$$
their respective composition is defined by $AfSys(\prm{\eta}\circ\eta)I_1\arw{\prm{\alpha}\boxdot\alpha}I_3\prm{\Phi}\Phi=
AfSys\prm{\eta}AfSys\eta I_1\arw{AfSys\prm{\eta}\alpha}AfSys\prm{\eta}I_2\Phi\arw{\prm{\alpha}\Phi}I_3\prm{\Phi}$ (cf. Definition~\ref{defn:7} and Lemma~\ref{lem:1}).

\par
Following the remarks after Definition~\ref{defn:10}, the direction of the natural transformation $\alpha$ in affine institution morphisms (cf. Definition~\ref{defn:12}) is changed.

\begin{defn}
 \label{defn:14.1}
 Given an affine theory $T$, \afinst{T} is the subcategory of \cat{AfInst} consisting of affine institutions $(\cat{S},T,I)$ (shortened to $(\cat{S},I)$) and their respective morphisms $(\Phi,1_T,\alpha)$ (shortened to $(\Phi,\alpha)$).~\sqed
\end{defn}

\begin{exmp}
 \label{exmp:5}
 \hfill\par
 \begin{enumerate}[(1)]
  \item If $\cat{B}=\cat{Frm}$, then the category \afinst{\mcal{P}_{\msf{2}}} provides a modification of the category
        \cat{TpInst} of~\cite{Sernadas1995}.
  \item The affine theory $\set\arw{|\mcal{P}|}\opm{\set}:=\set\arw{\mcal{P}}\opm{\cat{CBAlg}}\arw{\opm{|-|}}\opm{\set}$, in
        which $|-|$ is the obvious forgetful functor, gives the category \afinst{|\mcal{P}|}, which is a modification of the category \cat{ElInst} (recall the remark at the end of Definition~\ref{defn:9}).
        \sqed
 \end{enumerate}
\end{exmp}


\subsection{Affine institution spatialization and localification procedures}


\par
In this subsection, a possible approach to spatialization and localification procedures for affine institutions is presented. Following the ideas of~\cite{Sernadas1994,Sernadas1995}, we begin with some preliminary definitions.

\begin{defn}
 \label{defn:15}
 Let $T$ be an affine theory. A \emph{spatial affine $T$-institution} consists of:
 \begin{itemize}
  \item a category \cat{S} (of abstract signatures),
  \item a functor $\cat{S}\arw{I}\topt{T}$.
 \end{itemize}

 \par
 A \emph{spatial affine $T$-institution morphism} $(\cat{S}_1,I_1)\arw{(\Phi,\alpha)}(\cat{S}_2,I_2)$ consists of:
 \begin{itemize}
  \item a functor $\cat{S}_1\arw{\Phi}\cat{S}_2$,
  \item a natural transformation $I_1\arw{\alpha}I_2\Phi$.
 \end{itemize}

 \par
 \safinst{T} is the category of spatial affine $T$-institutions and their morphisms.
 \sqed
\end{defn}

\begin{defn}
 \label{defn:16}
 Let $\cat{X}\arw{T}\opm{\cat{B}}$ be an affine theory. A \emph{localic affine $T$-institution} consists of:
 \begin{itemize}
  \item a category \cat{S} (of abstract signatures),
  \item a functor $\cat{S}\arw{I}\opm{\cat{B}}$.
 \end{itemize}

 \par
 A \emph{localic affine $T$-institution morphism} $(\cat{S}_1,I_1)\arw{(\Phi,\alpha)}(\cat{S}_2,I_2)$ consists of:
 \begin{itemize}
  \item a functor $\cat{S}_1\arw{\Phi}\cat{S}_2$,
  \item a natural transformation $I_1\arw{\alpha}I_2\Phi$.
 \end{itemize}

 \par
 \lafinst{T} is the category of localic affine $T$-institutions and their morphisms.
 \sqed
\end{defn}

\begin{exmp}
 \label{exmp:6}
 If $\cat{B}=\cat{Frm}$, then the categories \safinst{\mcal{P}_\msf{2}} and \lafinst{\mcal{P}_\msf{2}} provide (slightly) modified versions of the categories of spatial and localic topological institution of~\cite{Sernadas1994,Sernadas1995}, respectively.
 \sqed
\end{exmp}

\par
From Theorems~\ref{thm:2} and \ref{thm:4} (employing their respective notations), the following results are obtained.

\begin{thm}
 \label{thm:7}
 \hfill\par
 \begin{enumerate}[(1)]
  \item There exists a full embedding \incl{\safinst{T}}{IE}{\afinst{T}} given by
         $IE((\cat{S}_1,I_1)\arw{(\Phi,\alpha)}(\cat{S}_2,I_2))=(\cat{S}_1,EI_1)\arw{(\Phi,E\alpha)}(\cat{S}_2,EI_2)$.
  \item $IE$ has a right-adjoint-left-inverse $\afinst{T}\arw{ISpat}\safinst{T}$,
        $ISpat((\cat{S}_1,I_1)\arw{(\Phi,\alpha)}(\cat{S}_2,I_2))=(\cat{S}_1,SpatI_1)\arw{(\Phi,Spat\alpha)}(\cat{S}_2,SpatI_2)$.
  \item \safinst{T} is isomorphic to a full coreflective subcategory of \afinst{T}.
 \end{enumerate}
\end{thm}
\begin{pf}
 For convenience of the reader, a sketch of the proof is given. Theorem~\ref{thm:2} supplies us with an adjoint situation $(\eta,\varepsilon):E\dashv Spat:\topsyst{T}\arw{}\topt{T}$, in which $1_{\topt{T}}\arw{\eta}SpatE$ and $ESpat\arw{\varepsilon}1_{\topsyst{T}}$ are natural transformations such that $Spat\varepsilon\circ\eta Spat=1_{Spat}$ and $\varepsilon E\circ E\eta=1_{E}$.

 \par
 Starting from $1_{\topt{T}}\arw{\eta}SpatE$, one obtains a natural transformation $1_{\safinst{T}}\arw{\Theta}ISpatIE$, which is given for every spatial affine $T$-institution $(\cat{S},I)$ by the following diagram
 $$
    \xymatrix{\cat{S} \ar[d]_-{I} \ar[rr]^-{1_{\cat{S}}} & & \cat{S} \ar[d]^{I}\\
              \topt{T} \ar@{=>}"1,3"^-{\eta I}  & \ar[l]^-{Spat} \topsyst{T} & \ar[l]^-{E} \topt{T}.}
 $$

 To see that $\Theta$ is indeed a natural transformation, it needs to be shown that given a spatial affine $T$-institution morphism
 $(\cat{S}_1,I_1)\arw{(\Phi,\alpha)}(\cat{S}_2,I_2)$, the diagram
 $$
   \xymatrix{(\cat{S}_1,I_1) \ar[d]_-{(\Phi,\alpha)} \ar[rr]^-{\Theta_{(\cat{S}_1,I_1)}} & & ISpatIE(\cat{S}_1,I_1)
              \ar[d]^{ISpatIE(\Phi,\alpha)}\\
             (\cat{S}_2,I_2) \ar[rr]_-{\Theta_{(\cat{S}_2,I_2)}} & & ISpatIE(\cat{S}_2,I_2)}
 $$
 commutes. However, the commutivity of the previous diagram follows because the next diagram
 $$
   \xymatrix{I_1 \ar[d]_-{\alpha} \ar[rr]^-{\eta I_1} & & SpatEI_1 \ar[d]^{SpatE\alpha}\\
             I_2\Phi \ar[rr]_-{\eta I_2\Phi} & & SpatEI_2\Phi,}
 $$
 is clearly commutative ($\eta$ is a natural transformation).

 \par
 In a similar way, one gets the natural transformation $IEISpat\arw{\Upsilon}1_{\afinst{T}}$. Moreover, to show commutativity of the triangle
 $$
   \xymatrix{ISpat \ar"2,3"_-{1_{ISpat}} \ar[rr]^-{\Theta ISpat} & & ISpatIEISpat \ar[d]^-{ISpat\Upsilon}\\
             & &ISpat,}
 $$
 one notices that $ISpat\Upsilon\circ\Theta ISpat=1_{ISpat}$ is nothing else than $Spat\varepsilon\circ\eta Spat=1_{Spat}$. Analogously, one arrives at $\Upsilon IE\circ IE\Theta=1_{IE}$, which concludes the proof.
 \qed
\end{pf}

\begin{thm}
 \label{thm:8}
 Let $T$ be an affine theory such that there exists an adjoint situation $(\eta,\varepsilon):T\dashv Pt:\loc{B}\arw{}\cat{X}$.
 \begin{enumerate}[(1)]
  \item There exists a functor $\afinst{T}\arw{ILoc}\lafinst{T}$, which is defined by
        $ILoc((\cat{S}_1,I_1)\arw{(\Phi,\alpha)}(\cat{S}_2,I_2))=(\cat{S}_1,LocI_1)\arw{(\Phi,Loc\alpha)}(\cat{S}_2,LocI_2)$.
  \item There exists a full embedding \incl{\lafinst{T}}{IE}{\afinst{T}} given by
        $IE((\cat{S}_1,I_1)\arw{(\Phi,\alpha)}(\cat{S}_2,I_2))=(\cat{S}_1,EI_1)\arw{(\Phi,E\alpha)}(\cat{S}_2,EI_2)$ such that $ILoc$ is a left-adjoint-left-inverse to $IE$.
  \item \lafinst{T} is isomorphic to a full reflective subcategory of \afinst{T}.
 \end{enumerate}
\end{thm}
\begin{pf}
 The proof is similar to that of Theorem~\ref{thm:7}, relying on Theorem~\ref{thm:4} instead of Theorem~\ref{thm:2}.
 \qed
\end{pf}

\par
Theorem~\ref{thm:7} gives an answer to the question regarding a spatialization construction for topological institutions posed in~\cite[p.~424]{Sernadas1995} (the authors though have refrained from establishing the respective coreflection).


\section{Conclusion and future work}


\par
Following the concept of topological institution of~\cite{Sernadas1994,Sernadas1995}, we introduced the notion of affine institution and showed its respective spatialization and localification procedures. Affine institutions seem to provide a better framework for elementary institutions and topological institutions because they do not require that the employed algebraic structures have to be frames. More precisely, while the authors of~\cite{Sernadas1994,Sernadas1995} impose the frame structure on the set of theories (certain ``closed" subsets of the set of sentences) of a given signature, which results in certain technical difficulties, we suggest the use of an arbitrary algebraic structure, which could be determined in each concrete case. Additionally, it is important to note that the idea of topological system, expressed in terms of affine systems, provides a convenient framework for, and is more general than, elementary institutions. This conclusion is in striking difference with the result of J.~T.~Denniston, A.~Melton and  S.~E.~Rodabaugh~\cite{Denniston2014}, which eventually represents lattice-valued topological systems as special lattice-valued institutions.

\par
We end this paper with four open problems.


\subsection{Institutions based in categories}


\par
This paper considered a simplification of the concept of institution, which relies on the category \set\ of sets instead of the category \cat{Cat} of categories. To incorporate the case of general institutions and their morphisms of Definitions~\ref{defn:9}, \ref{defn:10}, one could introduce the following concepts.

\begin{defn}
 \label{defn:17}
 Given a functor $\cat{Cat}\arw{T}\loc{Cat}$, \syst{T} is the comma category $(T\downarrow 1_{\loc{Cat}})$, concrete over the product category $\cat{Cat}\times\loc{Cat}$, whose objects (\emph{generalized $T$-affine systems}) are triples $(\cat{X},K,\cat{B})$, which are made by \loc{Cat}-morphisms $T\cat{X}\arw{K}\cat{B}$; and whose morphisms (\emph{generalized $T$-affine system morphisms})
 $(\cat{X}_1,K_1,\cat{B}_1)\arw{(F,\Phi)}(\cat{X}_2,K_2,\cat{B}_2)$ are $\cat{Cat}\times\loc{Cat}$-morphisms $(\cat{X}_1,\cat{B}_1)\arw{(F,\Phi)}(\cat{X}_2,\cat{B}_2)$, making the diagram
 $$
   \xymatrix{T\cat{X}_1 \ar[d]_{K_1} \ar[r]^{TF} & T\cat{X}_2 \ar[d]^{K_2}\\
             \cat{B}_1 \ar[r]_{\Phi} & \cat{B}_2}
 $$
 commute.
 \sqed
\end{defn}

\par
For convenience, a functor $\cat{Cat}\arw{T}\loc{Cat}$ is called a \emph{generalized affine theory}.

\begin{defn}
 \label{defn:18}
 A \emph{generalized affine institution} consists of:
 \begin{itemize}
  \item a category \cat{S} (of abstract signatures),
  \item a generalized affine theory $T$,
  \item a functor $\cat{S}\arw{I}\syst{T}$.
        \sqed
 \end{itemize}
\end{defn}

\begin{defn}
 \label{defn:19}
 A \emph{generalized affine institution morphism} $(\cat{S}_1,T_1,I_1)\arw{(F,\Phi,\alpha)}(\cat{S}_2,T_2,I_2)$ consists of:
 \begin{itemize}
  \item a functor $\cat{S}_1\arw{F}\cat{S}_2$,
  \item a functor $\syst{T_1}\arw{\Phi}\syst{T_2}$,
  \item a natural transformation $I_2F\arw{\alpha}\Phi I_1$, or, more specifically,
        $$
          \xymatrix{\cat{S}_1 \ar[d]_{I_1}
                    \ar[rr]^{F} & & \cat{S}_2 \ar@{=>}"2,1"_-{\alpha}
                    \ar[d]^{I_2}\\
                    \syst{T_1} \ar[rr]_{\Phi} & & \syst{T_2}.}
        $$
        \sqed
 \end{itemize}
\end{defn}

\par
Definitions~\ref{defn:18}, \ref{defn:19} provide the category \cat{GAfInst} of generalized affine institutions and their morphisms. The remarks after Definition~\ref{defn:9} say that the category \cat{GAfInst} is a generalization of the category \cat{Inst}. The two next problems then arise immediately.

\begin{prob}
 \label{prob:1}
 What are the spatialization and localification procedures for the category \cat{GAfInst}?
 \sqed
\end{prob}

\begin{prob}
 \label{prob:2}
 What are the morphisms between two given generalized affine theories $T_1$, $T_2$, which could induce functors between the respective categories \syst{T_1}, \syst{T_2}, so that one could use them (the morphisms) in the definition of the category \cat{GAfInst} (cf. the definition of the category \cat{AfInst})?
 \sqed
\end{prob}


\subsection{Variable-basis spatialization and localification procedures for institutions}


\par
Theorems~\ref{thm:7}, \ref{thm:8} deal with the fixed-basis (in the sense of~\cite{Hohle1999a}) category \afinst{T} of affine institutions over a given affine theory $T$. One could be interested, however, in the whole variable-basis (in the sense of~\cite{Rodabaugh1999a}) category \cat{AfInst}.

\begin{prob}
 \label{prob:3}
 What are the spatialization and localification procedures for the category \cat{AfInst}?
 \sqed
\end{prob}

\par
Note that one still does not have a description of possible spatialization and localification procedures for affine spaces and systems, based in different affine theories. More precisely, Theorems~\ref{thm:2}, \ref{thm:4} are fixed-basis w.r.t. the employed affine theory $T$.


\subsection{Bornological systems and institutions}


\par
Motivated by the concept of lattice-valued (fixed-basis) bornological space of~\cite{Abel2011},~\cite{Paseka2014} introduced the notion of bornological vector system as a bornological analogue of topological systems of S.~Vickers (see, e.g.,~\cite{Hogbe-Nlend1977} for an introduction into the theory of bornological spaces and their related concepts). In particular, the authors of~\cite{Paseka2014} provided a spatialization procedure for bornological vector systems (see also~\cite{Paseka}). Thus, the last open problem of this paper is as follows.

\begin{prob}
 \label{prob:4}
 What are the relationships (if any) between bornological systems and institutions? What kind of institutions could arise from bornological systems?
 \sqed
\end{prob}

\par
We notice that the second part of Problem~\ref{prob:4} could start a whole new development in the theory of institutions, providing, e.g., bornological institutions (\emph{\`{a} la} topological institutions of~\cite{Sernadas1994,Sernadas1995}).

\par
The above-mentioned open problems will be addressed in our future papers on the topic of institutions.


\section*{Acknowledgements}


\par
The authors would like to express their sincere gratitude to the anonymous referees of this paper for their helpful remarks and suggestions on its possible improvements.

\par
This is a preprint of the paper published in ``Fuzzy Sets and Systems". The final authenticated version of the paper is available online at: https://www.sciencedirect.com/science/article/pii/S0165011415003668.


\end{document}